\documentclass[a4paper,11pt]{article}
\usepackage[latin1]{inputenc}
\usepackage[T1]{fontenc}
\usepackage{amsfonts,amssymb}
\usepackage{pstricks}
\usepackage{amsthm}
\usepackage[all]{xy}
\usepackage[pdftex]{graphicx}
\usepackage{graphicx}   
\newcommand{\n}[0]{{\textbf N}}
\newcommand{\A}[0]{{\textbf A}}
\newcommand{\w}[0]{{\mathbb W}}
\renewcommand{\c}[0]{{\mathcal C}}
\renewcommand{\k}[0]{{\textbf K}}

\newtheorem{theo}{Theorem}[section]
\newtheorem{defi}{Definition}[section]
\newtheorem{prop}[theo]{Proposition}

\newtheorem{cor}[theo]{Corollary}
\newtheorem{lem}[theo]{Lemma}

\newtheorem{rem}[theo]{Remark}
\newtheorem{ex}{Example}

\title{Jet schemes and generating sequences of divisorial valuations in dimension two}

\author{Hussein Mourtada}

\begin {document}

\maketitle  
\begin{abstract}
Using the theory of jet schemes, we give a new approach to the description of a 
minimal generating sequence of a divisorial valuations on $\A^2.$ For 
this purpose, we show how one can recover the approximate roots of an 
analytically irreducible plane curve from the equations of its jet schemes. 
As an application, for a given divisorial valuation $v$ centered at the origin 
of $\A^2,$ we construct an algebraic embedding  
$\A^2\hookrightarrow \A^N,N\geq 2$ such that $v$ is the trace of 
a monomial valuation on $\A^N.$ We explain how results in this direction 
give a constructive approach to a conjecture of Teissier on resolution of 
singularities by one toric morphism.

\end{abstract}

\footnote{{\textbf{2010 Mathematics Subject Classification.} 
13A80,14E15,14E18,14M25.
This research was partially supported by the ANR-12-JS01-0002-01 SUSI.}}

\section{Introduction}

Let $X=\A^d=Spec~R,$ where $R=\k[x_1,\ldots,x_d]$ is a polynomial ring over an 
algebraically closed field $\k.$ The arc space of $X,$ that we denote by 
$X_{\infty},$ is the scheme whose $\k-$rational points are 
$$X_{\infty}(\k)=Hom_{\k}(Spec~\k[[t]],X).$$
We have a natural truncation morphism $X_{\infty}\longrightarrow X,$ that we 
denote by $\Psi_0.$ 
For $p\in \mathbb{\n}$ and $Y=V(I)\subset X$ a subvariety defined by an ideal 
$I$, we consider the subset of arcs in $X_{\infty}$ that have an order of 
contact $p$ with $Y,$ this is 

$$Cont^p(Y)=\{\gamma \in X_\infty\mid ord_t \gamma^*(I)=p\},$$
where $\gamma^*:R\longrightarrow \k[[t]]$ is the $\k-$algebra homomorphism 
associated with $\gamma$ and 
$$ord_t \gamma^*(I)=\mbox{min}_{h \in I} \big\{ ord_t \gamma^*(h) \big\}.$$ 

With an irreducible component $\w$ of $Cont^p(Y),$ which is contained in the 
fibre $\Psi_0^{-1}(0)$ above the origin, we associate a valuation 
$v_{\w}:R\longrightarrow \n$ as follows:
$$ v_{\w}(h)=\mbox{min}_{\gamma \in \w}\{ord_t\gamma^*(h)\},$$
for $h\in R.$ It follows from \cite{ELM} (see also \cite{dFEI}, \cite{Re}, prop. 
3.7 (vii)), that $v_{\w}$ is 
a divisorial valuation centered at the origin $0\in X,$ and that all divisorial 
valuations centered at $0\in X,$ can be obtained by this way. 

We are interested in determining a generating sequence of such a valuation, this 
is a sequence of elements of $R$ which determines completely the valuation. Let 
us explain what is such a sequence:\\
For $\alpha \in \n$, let $$\mathcal{P}_{\alpha}=\{ h \in R\mid v_{\w}(h)\geq 
\alpha \}.$$

Following \cite{T3}, we define the $\k$-graded algebra 
$$gr_{v_{\w}}R=\bigoplus_{\alpha \in \n} 
\frac{\mathcal{P}_{\alpha}}{\mathcal{P}_{\alpha+1}}.$$
We call $in_{v_{\w}}$ the natural map $$in_{v_{\w}}: R\longrightarrow  
gr_{v_{\w}}R,~~h\mapsto h ~~\mbox{mod}~~\mathcal{P}_{v_{\w}(h)+1} .$$

\begin{defi}\cite{S}
A generating sequence of $v_{\w}$ is a set of elements of $R$ such that their 
image by $in_{v_{\w}}$
generates $gr_{v_{\w}}R$ as a $\k$-algebra. 
\end{defi}

In this article we will give a new way to determine generating sequence of 
$v_{\w}$  in dimension $2,$ i.e. when 
$d=2.$
Traditionally, there are three approaches to determine such a generating 
sequence: \\

1)By studying the relations in the semigroup $v_{\w}(R)$ \cite{T3}. The new 
developements of this
theory in higher dimensions  treat only valuations with maximal rational rank 
\cite{T1}, \cite{T2}, which does not include divisorial valuations.\\

2)By considering curvettes \cite {S}: let $\pi$  be the composition of the 
minimal sequence of blow ups that produces the divisor defining $v_{\w}.$ Let 
$G$ be its dual graph, then a curvette is a curve which is an image of a 
transversal arc to a rupture divisor of $G.$ If we choose the equation of a 
curvette for every rupture divisor, plus the variables of $R,$ we obtain a 
generating sequence of $v_{\w}.$ This approach has not been generalized to 
higher dimensions, and this seems to be a difficult mission.\\

3)Maclane's method \cite{Mc}(see also \cite{AM},\cite{FJ}): A generating 
sequence is obtained by induction using euclidean division. The generalizations 
of this method to higher dimensions (\cite{V1},\cite{HOS},\cite{Ma}) do not 
produce elements in $R,$ which is essential for our applications. See also 
\cite{CV} for a comparable approach.\\ 

Our approach is based on the definition of a divisorial valuation that we gave 
above in terms of arcs (and jet schemes). It 
will enable us to build a generating sequence from the equations of the 
subset $\w$ of the arc space which defines the divisorial valuation. The 
construction of a generating sequence passes by the extraction of the 
approximate roots of a plane branch from its jet schemes.\\  

One motivating application that we will present, and which remains true for a 
particular type of divisorial valuations in higher dimensions \cite{Mo4}, is the 
following: Given  a divisorial valuation $v$ centered at $0 \in \A^2,$ we will 
determine an embedding 
$e: \A^2\hookrightarrow \A^n,$ (where $n$ depends on $v$) and a toric proper 
birational morphism  $\mu: X_{\Sigma}\longrightarrow \A^n$ such that:

\begin{displaymath}
    \xymatrix{
        \tilde{\A^2} \ar[r] \ar[d] & X_{\Sigma} \ar[d]^{\mu} \\
        \A^2 \ar@{^{(}->}[r]^e     & \A^n }
\end{displaymath}
$\bullet$ $X_{\Sigma}$ is a smooth toric variety (i.e $\Sigma$ is a fan which 
is obtained by a regular subdivision 
of the positive quadrant $\mathbb{R}_+^n,$ this quadrant is the cone that  
defines $\A^n$ as a toric variety),\\
$\bullet$  the strict transform $\tilde{\A^2}$ of $\A^2$ by $\mu$ is smooth,\\
$\bullet$ a toric divisor $E'$ (associated with one of the edges of $\Sigma$ 
and which is determined by the values of 
the elements in a generating sequence) intersects $\tilde{\A^2}$ 
transversally along a divisor $E;$ note that the valuation associated with 
$E'$ is monomial and is given by the weight vector corresponding to $E',$\\  
$\bullet$ the valuation defined by the divisor $E$ is $v.$ \\

Our goal is to use such a construction to answer constructively the following 
conjecture of Teissier \cite{T2}:\\

 For a subvariety $Y\subset \A^n,$ there exists an embedding  
$\A^n\hookrightarrow \A^N,N\geq n$ such that the singularities of $Y$ can be 
resolved by a birational proper toric map $Z\longrightarrow \A^N.$ \\
 
A solution of this problem in the case of quasi-ordinary singularities is given 
in \cite{GP}. A related result was proved in \cite{Te2}, but the author starts with
a given resolution of singularities.\\
 
 For a given  singular subvariety $Y\subset \A^n,$ our idea is to extract a 
finite number of significant divisorial valuations $v_1,\ldots,v_r$ on $\A^n$ 
from the jet schemes of $Y$ (this is to compare with the Nash map 
\cite{I},\cite{ELM}), then to embed as above  $\A^n$ in a larger affine space 
$\A^N$ in such a way that all the valuations $v_1,\ldots,v_r$ can be seen as the 
traces of  monomial valuations on $\A^N.$ If $v_1,\ldots,v_r,$ are well 
chosen, this should guarantee Newton non-degeneracy (\cite{AGS},\cite{Te1}) of 
$Y\subset \A^N$ and hence would give the desired embedding. There remains the 
subtle matter of detecting the valuations $v_1,\ldots,v_r,$ (see 
\cite{Mo3},\cite{LMR} for simple examples) and finding the embedding described 
above for general divisorial valuations. In \cite{Mo4} we present a progress in 
this last problem.\\  

This idea corresponds to an approach of resolution 
of singularities by one toric morphism which is different from the one 
suggested in \cite{GT}. Indeed in \textit{loc.cit.}, this resolution of an irreducible 
plane curve $\mathcal{C}$ is constructed by considering the curve valuation 
$\nu_{\mathcal{C}},$ while the approach suggested by this article is to study 
the divisorial valuations which are associated to special components of the jet 
schemes. The two approaches lead to the same result for plane branches, but 
bifurcate in higher dimensions.    

One application of the result of this article would be a resolution of 
singularities of a reducible plane curve with one toric morphism. This will be 
treated elsewhere.\\

I have found inspiration for this article in \cite{T2}, I am thankful to Bernard 
Teissier for all the explanations he gave to me about it, and for  
several corrections and suggestions he made about an earlier version of this 
article. I would also like to thank Pedro Gonz\'alez P\'erez, Monique 
Lejeune-Jalabert,  Mohammad Moghaddam, Patrick Popescu-Pampu and Matteo Ruggiero 
for several discussions about this subject.   

The article assumes some knowledge of valuations and toric geometry. This can be 
found respectively in \cite{V2} and \cite{AGS}.    

\section{Jet schemes}

Let $\k$ be an algebraically closed field of arbitrary characteristic.
Let $X$ be a $\k$-algebraic variety and let $m \in \mathbb{N}$. The functor $F_m 
:\k-\mbox{Schemes} \longrightarrow \mbox{Sets}$
which with an affine scheme defined by a $\k-$algebra $A$ associates
\\ $$F_m(Spec(A))=Hom_{\k}(Spec A[t]/(t^{m+1}),X)$$
is representable by a $\k-$scheme $X_m$ (\cite{EM},\cite{I}). $X_m$ is the $m-$th jet scheme 
of $X$, and $F_m$ is  isomorphic to its functor of points.
So we have the following bijection

\begin{equation}\label{bij}
Hom_{\k}(Spec A,X_m)\simeq Hom_{\k}(Spec A[t]/(t^{m+1}),X).
\end{equation}
If $X=Spec R$ is affine, then $X_m=Spec R_m$ is also affine, and by taking $A= 
R_m$ in the bijection  (\ref{bij}),
we obtain a universal morphism $\Lambda^*:R\longrightarrow R_m[t]/(t^{m+1}),$ 
which is the morphism associated 
to the image of the identity $id \in Hom_{\k}(X_m,X_m)$ by the bijection 
(\ref{bij}). For example, if 
$X=Spec\k[x_0,x_1],$ and $f\in\k[x_0,x_1]$ then 
$$X_m=Spec\k[x^{(0)}_0,x^{(0)}_1,\ldots,x^{(m)}_0,x^{(m)}_1]=Spec R_m,$$ and

\begin{equation}\label{uf}
\Lambda^*(f)=F^{(0)}+F^{(1)}t+\cdots +F^{(m)}t^m,
\end{equation}
where $F^{(i)}$ is the coefficient of $t^i$ in the expansion of
\begin{equation}
f(x^{(0)}_0+x^{(1)}_0t+\cdots x^{(m)}_0t^m,x^{(0)}_1+x^{(1)}_1t+\cdots 
x^{(m)}_1t^m).
\end{equation}
 
Note that since we are interested in the ideal generated by the $F^{(i)}$'s, in characteristic $0,$ we 
can reconstruct them in such a way that they are obtained by a derivation 
process, see proposition $2.3$ in \cite{Mo1}.

For $m,p \in \mathbb{N}, m > p$, the truncation homomorphism $A[t]/(t^{m+1}) 
\longrightarrow A[t]/(t^{p+1})$ induces
a canonical projection $\pi_{m,p}: X_m \longrightarrow X_p.$ These morphisms 
clearly satisfy $\pi_{m,p}\circ \pi_{q,m}=\pi_{q,p}$
for $p<m<q$, and they are affine morphisms, so that they define a projective 
system whose limit is a scheme that we denote $X_{\infty};$ it is the arc 
space of $X$.
\\ Note that $X_0=X$. We denote the canonical projection 
$\pi_{m,0}:X_m\longrightarrow X_0$ by $\pi_{m},$ and $\Psi_m$  the canonical 
morphisms $X_{\infty}\longrightarrow X_m.$ \\

\section[Generating sequences from jet schemes]{Minimal generating sequences of 
a curve valuation from the equations of jet schemes}

In \cite{Mo1} and \cite{LMR}, we have used the approximate roots to study the 
geometry of the jet schemes of plane branches 
and to obtain toric resolutions of singularities of these curves.
In this section we show how one can obtain a minimal generating sequence of the 
valuation defined by a plane branch, i.e  a curve valuation, from the jet 
schemes of the branch. Note that the graph
that we have introduced in \cite{Mo1} is not sufficient to determine this 
generating sequence. The invariants of the jet schemes that 
we consider below are finer and are not determined by the topological type.\\
Let $\mathcal{C}$ be a plane branch  defined by an irreducible power series $f 
\in \k [[x_0,x_1]],$ where $\k$ is
an algebraically closed field. We assume that $x_0=0$ (resp. $x_1=0)$ is 
transversal (resp. tangent) to $\mathcal{C},$ this always can be achieved by  a 
linear change of variables.
Let $\bar{\beta_0}, \cdots ,\bar{\beta_{g}}$ be the minimal system of generators 
of the semigroup $\Gamma (\mathcal{C})$ of $\mathcal{C}.$
Let $e_0=\bar{\beta_0}$ (this is also the multplicity of $\mathcal{C}$ at the origin)  and $e_i=gcd(e_{i-1},\bar{\beta_i}), i\geq 1$ (where 
$gcd$ is the great common divisor). Since the sequence of positive integers
$$e_0>e_1> \cdots >e_i> \cdots $$ is strictly decreasing, there exists $g \in 
\mathbb{N}$, such that $e_g=1$.
We set $$n_i:=\frac{e_{i-1}}{e_i}, m_i:=\frac{\beta_i}{e_i}, i=1, \cdots ,g$$
and by convention, we set $\beta_{g+1}=+\infty$ and $n_{g+1}=1$.
\\
We have that \\
\begin{enumerate}
\item $e_i=gcd(\bar{\beta_0}, \cdots ,\bar{\beta_{i}}),0 \leq i\leq g,$\\
\item For $1\leq i\leq g$, there exists a unique system of nonnegative integers 
$b_{ij},$
$0\leq j<i$ such that for $1\leq j<i$, $b_{ij}<n_j$  and $n_i 
\bar{\beta_{i}}=\Sigma_{0\leq j<i}b_{ij}\bar{\beta_{j}}$.
\end{enumerate}
To such a plane branch $\mathcal{C}=\{f=0\},$ we  associate a (curve) valuation
$$\nu_{\mathcal{C}}:\k[[x_0,x_1]]\longrightarrow \mathbb{N}\cup \infty,$$
which is positive on the maximal ideal $(x_0,x_1)$, by using local intersection 
multiplicity: $$\nu_{\mathcal{C}}(h)=\mbox{dim}\frac{\k[[x_0,x_1]]}{(f,h)},$$
for every $h \in \k[[x_0,x_1]].$ Note that $tr.deg (\nu_{\mathcal{C}})=0$ and 
$rank(\nu_{\mathcal{C}})=2$ (see \cite{FJ} page 17).\\

For an irreducible $h \in \k[[x_0,x_1]],$ we have that $h$ is, up to 
multiplication by a constant, of the form 
\begin{equation}\label{np}
h=(x_1^{n_h}-\alpha_h x_0^{m_h})^{\delta_h}+ \sum_{(a,b)} c_{ab}x_0^ax_1^b,
\end{equation}
where $m_h$ and $n_h$ are coprime, $\alpha_h \in \k^*, c_{ab} \in \k,$ and 
the points $(a,b)$ are strictly above the Newton polygon of $h$(\cite{CA}).

\begin{lem}\label{l1} Given $f,h$ in the form (\ref{np}) above, we have  
$x_1^{n_f}-\alpha_f x_0^{m_f}\not=x_1^{n_h}-\alpha_h x_0^{m_h}$ if and only if
$$\nu_{\mathcal{C}}(h)= 
\mbox{min}(\bar{\beta_0}m_h\delta_h,\bar{\beta_1}n_h\delta_h).$$
Moreover, we have that
$\left\{
\begin{array}{l}
  in_{\nu_{\mathcal{C}}}h=x_0^{m_h\delta_h}~~ \mbox{or}~~ 
x_1^{n_h\delta_h}~~\mbox{if}~~(m_f,n_f)\not=(m_h,n_h). \\
 in_{\nu_{\mathcal{C}}}h=(x_1^{n_h}-\alpha_h 
x_0^{m_h})^{\delta_h}~~\mbox{if}~~(m_f,n_f)=(m_h,n_h) ~~\mbox{and}~~\alpha_f 
\not= \alpha_h.
\end{array}
\right.$

\end{lem}
\begin{proof}
This follows from the the classical formula of the local intersection 
multiplicity: 
 $$\nu_{\mathcal{C}}(h)=ord_t h(x(t),y(t)),$$
 where $(x(t),y(t))$ is a special parametrization of $\mathcal{C}$ obtained by 
the Newton-Puiseux theorem \cite{CA}.
\end{proof}

Following \cite{Mo1}, we describe the irreducible components of the schemes of 
jets centered at $0,$ i.e. $\mathcal{C}_m^0:=\pi_m^{-1}(0),$ where 
$\pi_m:\mathcal{C}_m \longrightarrow \mathcal{C}$ is the canonical morphism.  \\
We set  $$Cont^e(x_0)_{m}(resp.Cont^{>e}(x_0)_{m}):=\{ \gamma \in \mathcal{C}_m 
\mid ord_tx_0 \circ \gamma=e(resp.>e)\},$$
then we can state :

\begin{theo}\label{a1}(Th. 4.9, \cite{Mo1})
Let $\mathcal{C}$ be a plane branch with $g$ Puiseux exponents. Let $m \in 
\mathbb{N}$. For $1\leq m < 
n_1\bar{\beta_1}+e_1$, $\mathcal{C}_m^0=Cont^{>0}(x_0)_m$ is irreducible.
For $q=\left[\frac{m-e_1}{n_1\bar{\beta_1}}\right]\geq 1,$ 
 the irreducible components of $\mathcal{C}_m^0$ are :
$$C_{m\kappa I}=\overline{Cont^{\kappa\bar{\beta_0}}(x_0)_m}$$ 
for $1\leq \kappa$ and $\kappa\bar{\beta_0}\bar{\beta_1}+e_1 \leq m,$\\
 $$C^j_{m\kappa v}=\overline{Cont^{\frac{\kappa\bar{\beta_0}}{n_j \cdots 
n_g}}(x_0)_m} $$ 
for $j=2, \cdots ,g,$  $1\leq \kappa$ and $\kappa \not \equiv 0~ mod ~n_j$ and 
such that
$\kappa n_1 \cdots n_{j-1}\bar{\beta_1}+e_1\leq m< \kappa\bar{\beta_j}$,\\
 
$$B_m= Cont^{>n_1q}(x_0)_m.$$
\end{theo}

We are interested in the following inverse system of irreducible components:

$$\cdots\longrightarrow C_{(\bar{\beta_0}\bar{\beta_1}+e_1+2)1 I}\longrightarrow  C_{(\bar{\beta_0}\bar{\beta_1}+e_1+1)1 I} 
\longrightarrow C_{(\bar{\beta_0}\bar{\beta_1}+e_1)1 I} \longrightarrow 
B_{\bar{\beta_0}\bar{\beta_1}+e_1-1}\longrightarrow 
B_{\bar{\beta_0}\bar{\beta_1}}.~~~~(\star)$$

Let $C_{m}:=\overline{Cont^{\bar{\beta_0}}(x_0)_m}$ (the notation $C_m$ will be 
used all over the paper). Let $\gamma_m$ be the generic point of $C_m.$
From corollary $4.2$ in \cite{Mo1}, we can see that for $m$ large enough, 
$$\mbox{ord}_t x_1 \circ\gamma_m(t)=\bar{\beta_1}.$$ 
Note that the only data we need to detect the inverse system $(\star)$ is 
the multiplicity $\bar{\beta_0}$ of the curve.
Indeed, the components in the system $(\star)$ are given by the closure of 
$Cont^{\bar{\beta_0}}(x_0)_m,$ for $m\geq \bar{\beta_0}\bar{\beta_1}-1.$ \\

 In the following lemma, we compute the intersection multiplicity  of two 
curves in terms of ideals of jet schemes. Our first goal is to give a new way 
to determine the initial part of an element $h\in \k[[x,y]],$ with respect to 
the valuation  $\nu_{\mathcal{C}}.$ This is achieved in corollary \ref{inc}.

 Let $D^m(x_0^{(\bar{\beta_0})})$ be the open subscheme of $\A^2_m$ defined by 
$x_0^{(\bar{\beta_0})}\not=0.$ Let $I_m$ be the ideal defining 
 $Cont^{\bar{\beta_0}}(x_0)_m$ in $D^m(x_0^{(\bar{\beta_0})})$ and let  $I_m^r$ 
be its radical. Let $h\in\k[[x,y]]$ be irreducible and $H^{(i)}$ be the 
coefficient of $t^i$ in $\Lambda^*(h)$ (see equation (\ref{uf})).\\
 \begin{rem} In what follows, unless stated otherwise, when we use the symbol 
$\equiv,$ we just
 want to replace elements which are congruent to zero by zero.
 
 \end{rem}
 
 \begin{lem}  \label{in} 
$\nu_{\mathcal{C}}(h)=l$ if and only if for $m>>0,$ we have that $H^{(i)}\equiv 
0$ mod $I_m^r,$ if $i<l,$
and  $H^{(l)} \not\equiv 0$ mod $I_m^r.$ 
 
 \end{lem}
 \begin{proof}
 If $\nu_{\mathcal{C}}(h)=l.$ We have that 
$\nu_{\mathcal{C}}(h)=ord_th(x_0(t),x_1(t)),$ for any good parametrization (i.e. a general point of
the curve corresponds to just one value of the parameter)
 $(x_0(t),x_1(t))$ of $\mathcal{C}.$ Let 
$x_0^{(0)},\ldots,x_0^{(i_m)},x_1^{(0)},\ldots,x_1^{(j_m)}$
 be the variables that intervene in the generators of $I_m^r.$ Note that 
$i_m,j_m<m.$ By definition of $I_m^r,$
 for any closed point 
$(a_0^{(0)},\ldots,a_0^{(i_m)},a_1^{(0)},\ldots,a_1^{(j_m)}) \in V(I_m^r)\subset 
Spec\k[x_0^{(0)},\ldots,x_0^{(i_m)},x_1^{(0)},\ldots,x_1^{(j_m)}],$
 there is a good parametrization of $\mathcal{C}$ of the form
 $$ 
(a_0^{(0)}+a_0^{(1)}t+\cdots+a_0^{(i_m)}t^{i_m}+\cdots,a_1^{(0)}+a_1^{(1)}
t+\cdots+a_1^{(j_m)}t^{j_m}+\cdots). $$
 
 It follows that 
 
$$ord_th(a_0^{(0)}+a_0^{(1)}t+\cdots+a_0^{(i_m)}t^{i_m}+\cdots,a_1^{(0)}+a_1^{
(1)}t+\cdots+a_1^{(j_m)}t^{j_m}+\cdots)=l,$$
 and so $H^{(i)}(a_0^{(0)},\ldots,a_0^{(i_m)},a_1^{(0)},\ldots,a_1^{(j_m)})=0$ 
for every $i<l,$ and 
 $$H^{(l)}(a_0^{(0)},\ldots,a_0^{(i_m)},a_1^{(0)},\ldots,a_1^{(j_m)})\not=0.$$ 
Hence, $H^{(i)}\equiv 0$ mod $I_m^r,$ for every $i<l,$
and  $H^{(l)} \not\equiv 0$ mod $I_m^r.$\\
The converse is straightforward. 
\end{proof}

\begin{rem}\label{est}
\begin{enumerate}
\item In the proof of lemma \ref{in}, the fact that for a closed point of 
$V(I_m^r) \subset 
Spec\k[x_0^{(0)},\ldots,x_0^{(i_m)},x_1^{(0)},\ldots,x_1^{(j_m)}],$ we find an 
arc which ``lifts'' this point, is not equivalent to
saying that any $m-$jet in the irreducible component defined by $I_m^r$ is 
liftable (which is not true). The reason is that we need more coordinates to 
define an $m-$jet, namely there remains to specify 
$x_0^{(i_m+1)},\ldots,x_0^{(m)},x_1^{(j_m+1)},\ldots,x_1^{(m)},$ which can be 
chosen freely, but for such a jet to be liftable, these coordinates should 
satisfy more equations.

\item  We can estimate the minimum $m$ that verifies lemma \ref{in}, by 
determining the variables that appears in the equations of jet schemes. We find
$$m=\kappa_h:=[l\frac{mult(f)}{mult(h)}],$$

where $mult$ denotes multiplicity, $l=\nu_{\mathcal{C}}(h)$ and the brackets 
$[~~]$ denote the integral part.
\end{enumerate}
\end{rem}

We continue with the settings of lemma \ref{in}. Since $H^{(l)} \not\equiv 0$ 
mod $I_{\kappa_h}^r,$ let $P\in \k[[x_0,x_1]],$ be the minimal part of $h$  such 
that 
$$(h-P)^{(i)}\equiv0~~\mbox{mod}~~ 
I_{\kappa_h}^r~~\mbox{for}~~ i\leq l.$$

This means that the terms (a term is a constant times a monomial in $x_0$ and $x_1$) of $P$ are terms of $h,$ and $P$ has the least number of terms with the above property. 
We thus obtain the following important corollary of lemma \ref{in}:
 
\begin{cor}\label{inc}

 We have that $$in_{\nu_{\c}}P=in_{\nu_{\c}}h.$$ Moreover, $P$ is the minimal part of $h$ achieving this equality.

\end{cor}
\begin{proof}
 It follows from the definition of $P$ and from lemma \ref{in} that 
$\nu_{\c}(h-P)>\nu_{\c}(h),$ and the assertion follows.\\
\end{proof}
\begin{rem}
 
\end{rem}

\begin{ex}We assume the characteristic of $\k$ is zero, which makes the 
computation easier.
\begin{enumerate}
\item Let $\mathcal{C}=\{f=(x_1^2-x_0^3)^2-x_0^6x_1=0\},$ and let 
$h=(x_1^2-x_0^3)^2-4x_0^5x_1-x_0^7.$
We have that $\nu_{\c}(h)=26.$ We can see this by applying lemma \ref{in}, 
indeed:
$$I^r_{26}=(x_0^{(0)},\ldots,x_0^{(3)},x_1^{(0)},\ldots,x_1^{(5)},{x_1^{(6)}}^2-
{x_0^{(4)}}^3,
2x_1^{(6)}x_1^{(7)}-3{x_0^{(4)}}^2x_0^{(5)})\subset (R_m)_{x_0^{(4)}},$$
where $(R_m)_{x_0^{(4)}}$ is the ring $R_m$ localized by $x_0^{(4)}.$ Note that, 
since in this example we have that $f$ and $h$ have the same multiplicity, we have 
$\kappa_h=\nu_{\c}(h).$
We observe that for every $i<26,~H_i\equiv 0$ modulo $I^r_{26},$ and 
$$H^{(26)}\equiv -4{x_0^{(4)}}^5x_1^{(6)}\not\equiv 0~~\mbox{mod}~~ I^r_{26}.$$
From corollary \ref{inc}, we deduce that $in_{\nu_{\c}}h=-4x_0^5x_1.$ 

\item Let $\mathcal{C}=\{f=(x_1^2-x_0^3)^2-x_0^6x_1\},$ and let $h=x_1^2-x_0^3.$
We have that $\nu_{\c}(h)=15,$ and therefore $\kappa_h=30.$ We have that
$$I^r_{30}=(x_0^{(0)},\ldots,x_0^{(3)},x_1^{(0)},\ldots,x_1^{(5)},H^{(12)}, 
H^{(13)},H^{(14)},
{H^{(15)}}^2-{x_0^{(4)}}^6x_1^{(6)})\subset (R_m)_{x_0^{(4)}},$$

where $H^{(12)}={x_1^{(6)}}^2-{x_0^{(4)}}^3$ and 
$H^{(13)}=2x_1^{(6)}x_1^{(7)}-3{x_0^{(4)}}^2x_0^{(5)}.$
We observe that for every $i<15,~H^{(i)}\equiv 0$ modulo $I^r_{30},$ and 
$$H^{(15)}\not\equiv 0~~\mbox{mod}~~ I^r_{30}.$$
From corollary \ref{inc}, we deduce that $in_{\nu_{\c}}h=h.$

\end{enumerate}
\end{ex}

Let us have a look at the equations of jet schemes. It follows from corollary 
$4.2$ in \cite{Mo1} that 
\begin{equation}\label{mon}I_{\bar{\beta_0}\bar{\beta_1}-1}=(x_0^{(0)},\ldots,
x_0^{(\bar{\beta_0}-1)},x_1^{(0)},\ldots,x_1^{(\bar{\beta_1}-1)}).
\end{equation}
We get from the same corollary that 

\begin{equation}\label{mon1}F^{(\bar{\beta_0}\bar{\beta_1})}\equiv 
({x_1^{(\bar{\beta_1})}}^{n_1}-c{x_0^{(\bar{\beta_0})}}^{m_1})^{e_1}~\mbox{mod}
~I_{\bar{\beta_0}\bar{\beta_1}-1},
\end{equation}
for some $c\in \k,c\not=0.$\\

\begin{rem}\label{n1}
Note that the equations (\ref{mon}) and (\ref{mon1}) are conditional on the 
hypothesis we have made on the variables
$x_0$ and $x_1.$ These variables  permit the best approximation of the valuation 
$\nu_{\mathcal{C}}$  by a monomial valuation, namely the monomial valuation 
$\nu_1$ which is determined by $\nu_1(x)=\nu_{\mathcal{C}}(x)$ and 
$\nu_1(y)=\nu_{\mathcal{C}}(y)$. Note that if we begin with any choice of 
variables, we can use jet schemes to detect variables verifying this property.
\end{rem}

We now give the steps of an algorithm that determine the minimal generating 
sequence. This will be guided by the fact that we can detect the initial part of
a function with respect to $\nu_\mathcal{C}$ from the equations of the jet schemes 
of $\mathcal{C};$ this follows from lemma \ref{in} and corollary \ref{inc}. So we will determine
algorithmically  elements in $\k[x,y]$ whose images by the universal morphism $\Lambda^*$
(see equation  \ref{uf}) generate the equations of the families of jets that define the 
valuations $\nu_\mathcal{C}.$\\

If $e_1=1,$ then a minimal generating sequence of $\nu_\mathcal{C}$ is given by 
$x_0,x_1$ and $f$ itself. We assume that $e_1>1.$ 

We set $x_{2,0}=x_1^{n_1}-x_0^{m_1}$ and to every 
$C_{m}:=\overline{Cont^{\bar{\beta_0}}(x_0)_m}$ in $(\star),$ we assign a vector 
$v^{3,0}_m=v^{3,0}(C_{m}) \in  \mathbb{N}^3$  as follows:

$$v^{3,0}_m=(\mbox{ord}_t x_0 \circ \gamma_m(t),\mbox{ord}_t x_1 
\circ\gamma_m(t),\mbox{ord}_t x_{2,0} \circ\gamma_m(t)), $$
where $\gamma_m$ is the generic point of $C_{m}.$ Let 
$$\mu_{2,0}=\mbox{min}\{m\geq \bar{\beta_0}\bar{\beta_1} \mid 
\mbox{codim}(C_{m+1}) > \mbox{codim}(C_m)$$
$$\mbox{and} ~~~~v^{3,0}_m=v^{3,0}_{m+1}\}.$$

Let $$F^{(\mu_{2,0}+1)} \equiv ~Q^{l} ~~\mbox{mod}~~I^r_{\mu_{2,0}},$$
for some reduced polynomial $Q$ and positive integer $l;$ note that $Q$ is a 
non zero polynomial because the equation $F^{(\mu_{2,0}+1)}$ forces the 
inequality
$\mbox{codim}(C_{\mu_{2,0}+1}) > \mbox{codim}(C_{\mu_{2,0}})).$\\

We then have two cases: \\
\underline{Case 1:} $l=e_1.$ \\

\textbf{Claim 1}\\

If $l=e_1,$ then we have that $$Q-x_{2,0}^{(\frac{\mu_{2,0}+1}{e_1})} \equiv 
Q'~~\mbox{mod}~~I^r_{\mu_{2,0}},$$  
where $Q'(x_0^{(\bar{\beta_0})},x_1^{(\bar{\beta_1})})$ is a polynomial in the 
variables $x_0^{(\bar{\beta_0})}$ and $x_1^{(\bar{\beta_1})}.$ 

We then define $$x_{2,1}=x_{2,0}+Q'(x_0,x_1),$$
and 
$$v^{3,1}_m=(\mbox{ord}_t x_0 \circ \gamma_m(t),\mbox{ord}_t x_1 
\circ\gamma_m(t),\mbox{ord}_t x_{2,1} \circ\gamma_m(t)).$$

\underline{Case 2:} If $l=l_2<e_1,$ and $l_2=1$ we stop. \\

\textbf{Claim 2}\\

If $1<l_2<e_1,$ then we have that 
$$Q-{x_{2,0}^{(\frac{\mu_{2,0}+1}{e_1})}}^{\frac{e_1}{l_2}} \equiv 
Q'~~\mbox{mod}~~I^r_{\mu_{2,0}},$$ 
where $Q'(x_0^{(\bar{\beta_0})},x_1^{(\bar{\beta_1})})$ is a polynomial in the 
variables $x_0^{(\bar{\beta_0})}$ and $x_1^{(\bar{\beta_1})}.$\\

We then set $x_2:=x_{2,0},\mu_{2}:=\mu_{2,0}$ and define 
$$x_{3,0}=x_2^{\frac{e_1}{l_2}}+Q'(x_0,x_1),$$
and 
$$v^{4,0}_m=(\mbox{ord}_t x_0 \circ \gamma_m(t),\mbox{ord}_t x_1 
\circ\gamma_m(t),\mbox{ord}_t x_2 \circ\gamma_m(t),
\mbox{ord}_t x_{3,0} \circ\gamma_m(t)).$$
We assume that we have recursively determined 
$(x_2,\ldots,x_{i-1},x_{i,j}),(e_1,l_2,\ldots,l_{i-1})$ and 
$(\mu_2,\ldots,\mu_{i-1},\mu_{i,j-1})$(if $j=0,$ we set 
$\mu_{i,j-1}=\mu_{i-1}).$\\
We define 
$$v^{i,j}_m=(\mbox{ord}_t x_0 \circ \gamma_m(t),\mbox{ord}_t x_1 
\circ\gamma_m(t),\ldots,\mbox{ord}_t x_{i,j} \circ\gamma_m(t)),$$   
   
and 
$$\mu_{i,j}=\mbox{min}\{m\geq \mu_{i,j-1}+1 \mid \mbox{codim}(C_{m+1}) > 
\mbox{codim}(C_m)$$
$$\mbox{and} ~~~~v^{i,j}_m=v^{i,j}_{m+1}\}.$$
Let $$F^{(\mu_{i,j}+1)} \equiv ~Q^l ~~\mbox{mod}~~I^r_{\mu_{i,j}},$$
for some reduced polynomial $Q$ and positive integer $l;$ note as above that 
$Q$ is a non zero polynomial because the equation $F^{(\mu_{i,j}+1)}=0$ forces 
the inequality $\mbox{codim}(C_{\mu_{i,j}+1}) > \mbox{codim}(C_{\mu_{i,j}}).$\\

We then  have two cases: \\
\underline{Case 1:} If $l=l_{i-1}.$\\

\textbf{Claim 1} continues\\

Then we have that $$Q-x_{i,j}^{(\frac{\mu_{i,j}+1}{l_{i-1}})} \equiv 
Q'~~I^r_{\mu_{i,j}},$$  
where $Q'$ is a polynomial in $x_0^{(\bar{\beta_0})},x_1^{(\bar{\beta_1})},
{x_{2}}^{(\frac{\mu_{2}+1}{e_1})},\ldots,{x_{i-1}}^{(\frac{\mu_{i-1}+1}{l_{i-2}}
)}.$\\
We then define $$x_{i,j+1}=x_{i,j}+Q'(x_0,x_1,\ldots,x_{i-1}),$$
and 
$$v^{i,j+1}_m=(\mbox{ord}_t x_0 \circ \gamma_m(t),\ldots,\mbox{ord}_t x_{i,j+1} 
\circ\gamma_m(t)).$$

\underline{Case 2:} If $l=l_{i}<l_{i-1},$ and $l_{i}=1$ we stop. If 
$1<l_{i}<l_{i-1}.$\\

\textbf{Claim 2} continues\\

Then we have that 
$$Q-{x_{i,j}^{(\frac{\mu_{i,j}+1}{e_1})}}^{\frac{l_{i-1}}{l_i}} \equiv 
Q'~~\mbox{mod}~~I^r_{\mu_{i,j}},$$ 
where $Q'$ is a polynomial in $x_0^{(\bar{\beta_0})},x_1^{(\bar{\beta_1})},
{x_{2}}^{(\frac{\mu_{2}+1}{e_1})},\ldots,{x_{i-1}}^{(\frac{\mu_{i-1}+1}{l_{i-2}}
)}.$\\

We then set $x_i:=x_{i,j},\mu_{i}:=\mu_{i,j}$ and define 
$$x_{i+1,0}=x_i^{\frac{l_{i-1}}{l_i}}+Q'(x_0,x_1,\ldots,x_{i-1}).$$

\begin{rem} If we want the elements of a generating sequence to be polynomials 
(which is more consistent with the terminology
key polynomials), then we might need an infinite number of elements to form a 
generating sequence (\cite{FJ}). These polynomials can be found by continuing 
the same algorithm, without stopping if we reach $l_g=1,$ but only if we reach 
$f,$ a case which occurs after finitely many steps if and only if $f$ is a 
polynomial. Here, we permit
as in \cite{T1}, elements in the ring $\k[[x_0,x_1]].$ Hence, even if $f$ is a 
power series and not a polynomial, we will take $f$
as an element of a ``minimal'' generating sequence. In that case, we can think 
$f$ as a limit key polynomial.   
\end{rem}

\textbf{Proof of claim 1}

By definition of $F^{(\mu_{2,0}+1)}$, the term $Q'$ comes from a polynomial $P$ 
such that the terms of $P^{e_1}$  appears in $f.$ More precisely,
$$Q'\equiv P^{(\frac{\mu_{2,0}+1}{e_1})}~~\mbox{mod}~~I^r_{\mu_{2,0}}.$$
By construction we then have that 
$$x_{2,0}^{(\frac{\mu_{2,0}+1}{e_1})}\equiv P^{(\frac{\mu_{2,0}+1}{e_1})} 
~~\mbox{mod}~~I^r_{\mu_{2,0}},$$
and both members are not congruent to $0$ modulo $I^r_{\mu_{2,0}+1}$ (because 
the codimension of the irreducible component of the 
$(\mu_{2,0}+1)-$jet scheme that we are considering increases).
We deduce from lemma \ref{in} that 
$$\nu_{\mathcal{C}}(x_{2,0}-P)>\nu_{\mathcal{C}}(x_{2,0})=\nu_{\mathcal{C}}
(P)=\frac{\mu_{2,0}+1}{e_1},$$
which implies that,
$in_{\nu_{\mathcal{C}}}  x_{2,0}=x_{2,0}=in_{\nu_{\mathcal{C}}}P.$ Indeed, any 
polynomial whose terms are also terms of $x_{2,0},$
namely $x_1^{n_1}$ and $x_0^{m_1},$ have value less than 
$\nu_{\mathcal{C}}(x_{2,0}).$

We have that $P$ is of the form 
$$P=P_1^{a_1}\cdots P_s^{a_s}.$$
This follows from the fact that the residue field of $\nu_{\mathcal{C}}$ is 
$\k,$ since $tr.deg(\nu_{\mathcal{C}})=0$ and $\k$ is algebraically closed.
We want to prove that the $in_{\nu_{\mathcal{C}}}P_j$'s are monomials in $x_0$ 
and $x_1$ for every $j.$ If not, then by lemma \ref{l1} we have 
that $$P_j=(x_1^{n_1}- \alpha_fx_0^{m_1})^{\delta_{P_j}}+ \sum 
c_{ab}x_0^ax_1^b$$
where $(a,b)$ is above the Newton polygon of $P_j.$ If $(x_1^{n_1}- 
\alpha_fx_0^{m_1})^{\delta_{P_j}}$ is a part of $in_{\nu_{\mathcal{C}}}P_j,$
this implies that $\nu_{\mathcal{C}}(P_j)\geq \nu_{\mathcal{C}}(x_{2,0)},$ and 
the equality follows from 
$in_{\nu_{\mathcal{C}}}  x_{2,0}=in_{\nu_{\mathcal{C}}}P.$ We deduce that 
$\delta_{P_j}=1.$ Then $in_{\nu_{\mathcal{C}}}P=in_{\nu_{\mathcal{C}}}P_j$ 
contains 
$x_{2,0},$  which contradicts the form of the equation (\ref{np}) for $f.$ It 
follows that $(x_1^{n_1}- cx_0^{m_1})^{\delta_{P_j}}$ is not  a part of 
$in_{\nu_{\mathcal{C}}}P_j,$ and we deduce by lemma $\ref{l1}$ that 
$in_{\nu_{\mathcal{C}}}P_j$ is a sum of monomials in $x_0$ and $x_1.$\\

Let's prove he remaining part of the claim 1; the proof is by induction on $i,$ 
we assume that the claim is true till $i-1:$ again  the term $Q'$ (in "claim 1 
continues") comes from a polynomial $P$ 
such that the terms of $P^{l}$  appears in $f.$ We have that $P$ of the form 
$$P=P_1^{a_1}\cdots P_s^{a_s},$$
where $P_r$ is irreducible for $r=1,\ldots,s.$ This again follows from the fact 
that $tr.deg(\nu_{\mathcal{C}})=0$. Note that as above 
$$\nu_{\mathcal{C}}(P)=\frac{\mu_{i,j}+1}{l},$$ and we will have 
$\nu_{\mathcal{C}}(P_r)\leq \frac{\mu_{i-1}}{l_{i-1}}.$
It follows from corollary \ref{inc} and from the hypothesis of induction that  
$in_{\nu_{\mathcal{C}}}P_r$ is a polynomial
in 
$x_0^{(\bar{\beta_0})},x_1^{(\bar{\beta_1})},{x_{2}}^{(\frac{\mu_{2}+1}{e_1})},
\ldots,{x_{i-1}}^{(\frac{\mu_{i-1}+1}{l_{i-2}})}.$
The proof of claim $2$ is similar to the proof of claim 1.

\begin{theo}\label{gc}
We have that 

\begin{enumerate}
\item For $i=2,\ldots,g,$ $\mu_i=e_{i-1}\bar{\beta_i}-1,$
and $l_i=e_i.$  Therefore $l_g=1$ and the algorithm stops at 
$\mu_{g}=e_{g-1}\bar{\beta_{g}}.$
\item $x_0, x_1,\ldots,x_g,f$ is a minimal generating sequence of 
$\nu_{\mathcal{C}}.$
\end{enumerate}
\end{theo}
\begin{proof}
The first part follows from  the formula for the codimension of $C_m$ in 
proposition 4.7 of \cite{Mo1} and the construction of the $\mu_i'$s.
We also recover that $\nu_{\mathcal{C}}(x_i)=\bar{\beta_i},i=0,\ldots,g.$    
The second part follows from corollary \ref{inc} and the description of the 
equations defining $C_m$ in terms of the equations of the jet schemes of the 
curves defined by $x_i,i=1,\ldots,g.$ Note that according to claim 1, 
$in_{\nu_{\mathcal{C}}}x_{i,j}$ is generated by $x_0,\ldots,x_i.$

\end{proof}

\begin{ex}
Let $f=((x_1^2-x_0^3-x_0^4)^2-x_0^8x_1)^2-x_0^3x_1(x_1^2-x_0^3-x_0^4),$ and  let 
$\mathcal{C}$ be the curve defined by $f.$
We have that $e_1=4,$ 
$x_{2,0}=x_1^2-x_0^3$ and 
$\mu_{2,0}=127.$ Let 
$$F^{(\mu_{2,0}+1)} \equiv ~Q^{l} ~~\mbox{mod}~~I^r_{\mu_{2,0}},$$ 
then $Q={x_{2,0}}^{(32)}-{{x_0}^{(8)}}^4$ and $l=4=e_1,$ hence we define 
$$x_{2,1}=x_{2,0}-x_0^4=x_1^2-x_0^3-x_0^4.$$
We have that $\mu_{2,1}=151.$ Let $$F^{(\mu_{2,1}+1)} \equiv ~Q^{l} 
~~\mbox{mod}~~I^r_{\mu_{2,1}},$$
then $Q={{x_{2,1}}^{(38)}}^2-{{x_0}^{(8)}}^8{x_1}^{(12)}$ and $l=l_2=2<e_1.$  
Since $l_2=2<e_1,$ we set $\mu_{2} :=\mu_{2,1},x_{2,1}=x_2,$ and we define $$ 
x_{3,0}= x_2^2-x_0^8x_1=(x_1^2-x_0^3-x_0^4)^2-x_0^8x_1.$$ \\
We have that $\mu_{3,0}=153,$ and we find that $l_3=1<l_2,$ hence we set  
$\mu_{3}:=\mu_{3,0},x_3=x_{3,0}$ and we stop.  
A minimal system of generator is then given by $x_0,x_1,x_2,x_3 $ and $f.$
\end{ex}

\section{Generating sequences of divisorial valuations}

We now apply the results of the previous section to determine from the jet 
schemes a minimal generating sequence for a divisorial valuation centered at 
the origin of $\A^2.$ The key point is that in 
dimension $2,$ a divisorial valuation $\nu_E$ which is determined by a divsor 
$E$ is defined by an irreducible component of $Cont^p(\mathcal{C}),$ where $p 
\in \mathbb{N}$ and 
$\mathcal{C}$ is an analytically irreducible plane curve. More precisely, the 
valuation is given by an irreducible component of $\mathcal{C}_{p-1}$ which is of 
the type  $C_{p-1}$ (see the definition of $C_m$ after theorem \ref{a1}) for 
$p\geq n_g\bar{\beta}_g,$ where $\bar{\beta}_0,\ldots,\bar{\beta}_g$ give a 
minimal system of generators of the semigroup $\Gamma(\mathcal{C}).$ Note that 
these numbers (the $\bar{\beta}_i $'s) are also extracted from the jet schemes, 
this is the first part of theorem \ref{gc}.  \\ 
The existence of $\mathcal{C}$ follows for instance from theorem 2.7 in 
\cite{LMR}:  $\mathcal{C}$ is chosen to be a curvette of $E.$ Recall that 
$\mathcal{C}$ is a curvette of $E,$ if there exists $\pi: X\longrightarrow 
\A^2,$ a composition of point blow ups above the origin, where $E$ is an 
irreducible component of the exceptional divisor of $\pi$ and the strict 
transform of $\mathcal{C}$ by $\pi$ is smooth and transversal to $E$ at a point 
which is not an intersection of $E$ with an other component of the exceptional 
divisor, i.e. a free point \cite{GB},\cite{FJ}.      

We will obtain a generating sequence of $\nu_E$ from the equations of the jet 
schemes of the curvette $\mathcal{C},$ more precisely 
from the irreducible component $C_{p-1}.$ There are two cases:\\

\underline{If $p=n_g\bar{\beta}_g,$} \\

let $x_2,\ldots, x_g$ be constructed by the algorithm of the previous section. 
Then a minimal generating sequence of the valuation $\nu_E$ is given by 
$x_0,\ldots,x_{g}.$ This follows from the definition of  $\nu_E$ in terms
of jet schemes. Indeed, $C_{p-1}$ gives rise to an irreducible component 
$\mathbb{W}$ of $Cont^p(\mathcal{C})$ (see the discussion after theorem 3.2 in 
\cite{Mo2}), and we have that  $$ v_{\w}(h)=\mbox{min}_{\gamma \in 
\w}\{ord_t\gamma^*(h)\},$$
for $h\in R=\k[x_0,x_1].$\\

\underline{If $p>n_g\bar{\beta}_g,$}\\

then we need to continue the algorithm in the previous section. Assume that we 
have constructed $x_0,\ldots,x_{g-1},$
hence we have  $$F^{(n_g\bar{\beta}_g)} \equiv ~Q 
~~\mbox{mod}~~I^r_{n_g\bar{\beta}_g-1},$$
for some reduced polynomial $Q$ (note that we do not have a power of $Q$ because 
we have reached the step where $l=1).$
We have that $$Q-{x_{g}^{(\bar{\beta}_g)}}^{n_g} \equiv 
Q'~~\mbox{mod}~~I^r_{n_g\bar{\beta}_g-1},$$ 
where $Q'$ is a polynomial in $x_0^{(\bar{\beta_0})},x_1^{(\bar{\beta_1})},
{x_{2}}^{(\bar{\beta}_2)},\ldots,x_{g-1}^{(\bar{\beta}_{g-1})}.$ We then define 
$$x_{g+1,0}=x_g^{n_g}+Q'.$$

We define 
$$v^{g+2,0}_m=(\mbox{ord}_t x_0 \circ \gamma_m(t),\mbox{ord}_t x_1 
\circ\gamma_m(t),\ldots,\mbox{ord}_t x_{g+1,0} \circ\gamma_m(t)),$$   
   
and 
$$\mu_{g+1,0}=\mbox{min}\{ n_g\bar{\beta}_g\leq m <p \mid \mbox{and} 
~~~~v^{g+2,0}_m=v^{g+2,0}_{m+1}\}.$$

We have not imposed any condition on the codimension in the definition of 
$\mu_{g+1,0}$ because for $m\geq n_g\bar{\beta}_g$ the codimension of $C_m$ 
grows by $1$ when $m$ grows by $1$ (proposition 4.7 in \cite{Mo1}).\\
 If $\mu_{g+1,0}=p-1,$ then a minimal generating sequence of $\nu_E$ is given by 
$$x_0,\ldots,x_{g+1}:=x_{g+1,0}.$$ If not, let $$F^{(\mu_{g+2,0}+1)} \equiv 
~Q~~\mbox{mod}~~I^r_{\mu_{g+2,0}},$$ for some reduced polynomial $Q.$ We have 
that 
 $$Q-x_{g+1,0}^{(\mu_{g+1,0}+1)} \equiv Q'~~\mbox{mod}~~I^r_{\mu_{g+1,0}},$$ 
where $Q'$ is a polynomial in $x_0^{(\bar{\beta_0})},x_1^{(\bar{\beta_1})},
{x_{2}}^{(\bar{\beta}_2)},\ldots,x_{g}^{(\bar{\beta}_{g-1})}.$ We then define 
$$x_{g+1,1}=x_{g+1,0}+Q'.$$  
Again, we define as above 
$v^{g+2,1}_m,\mu_{g+1,1},x_{g+1,2},\ldots,v^{g+2,j}_m,\mu_{g+1,j},$ until we 
have $\mu_{g+1,j}=p-1$ (note that $\mu_{g+1,i+1}>\mu_{g+1,i},i\geq0).$ Then a 
minimal generating sequence of $\nu_E$ is given by 
$$x_0,\ldots,x_{g+1}:=x_{g+1,j}.$$ 
Note that $\nu_E(x_g)=\bar{\beta}_g,\nu_E(x_{g+1})=p$ and all the 
$x_i$'s are polynomials in $\k[x_0,x_1].$ Actually if we let 
$\mathcal{D}=\{x_{g+1}=0\},$ then it follows from the 
definitions of $\nu_E$ and $\mathcal{D}$ that for an irreducible $h\in 
\k[[x,y]],$ we have that $$\nu_E(h)=\nu_{\mathcal{D}}(h)$$ and 
the initial part $in_{\nu_E}(h)=in_{\nu_{\mathcal{D}}}(h)$ is a polynomial in 
$x_0,\ldots,x_g,x_{g+1,j-1},$ unless  if
$in_{\nu_E}(h)=x_{g+1}^r$ is a power of $x_{g+1},$ in which case we have that 
$\nu_E(h)=rp.$ Note that $x_{g+1,j-1}$ is a polynomial in the variables 
$x_0,\ldots,x_{g+1}.$\\    

We now assume that for a divisorial valuation $\nu_E,$ defined by the 
irreducible component $C_{p-1}$ of the $(p-1)-th$ jet scheme of an irreducible 
curve $\mathcal{C},$ we have determined  $x_0,\ldots,x_g$ a minimal generating 
sequence as above. Then,  by construction, we have that for $i=2,\ldots,g,$ 
there exist polynomials $f_i$ such that
$$x_{i}=f_i(x_0,\ldots,x_{i-1}).$$
We will use this to prove the following proposition which is the goal of this 
article.

\begin{prop}\label{emb}
There exists an embedding $e: \A^2\hookrightarrow \A^{g+1},$  and a toric proper 
birational  morphism  $\mu: X_{\Sigma}\longrightarrow \A^{g+1}$ such that:

\begin{displaymath}
    \xymatrix{
        \tilde{\A^2} \ar[r] \ar[d]_{\eta} & X_{\Sigma} \ar[d]^{\mu} \\
        \A^2 \ar@{^{(}->}[r]^e     & \A^{g+1} }
\end{displaymath}
\begin{enumerate}
\item $X_{\Sigma}$ is smooth, i.e the fan $\Sigma$ is a regular subdivision of 
$\mathbb{R}_+^{g+1},$ and the vector 
$$v_{\nu_E}:=(\nu_E(x_0),\ldots,\nu_E(x_g))$$
 is an edge of a cone which belongs to $\Sigma,$
\item  the strict transform $\tilde{\A^2}$ of $\A^2$ by $\mu:X_{\Sigma} 
\longrightarrow  \A^{g+1} $ is smooth,
\item The divisor $E'\subset X_{\Sigma},$ which corresponds to the vector 
$v_{\nu_E},$ intersects $\tilde{\A^2}$ transversally along a divisor $E,$
\item the valuation defined by the divisor $E$ is $v.$ 
\end{enumerate}
\end{prop}

\begin{proof}
The functions $f_i$'s provide an embedding  $\A^2\hookrightarrow \A^{g+1},$ 
which is the geometric counterpart of the following morphism 
$$ \k[x_0,x_1,y_1,\ldots,y_{g}]\longrightarrow 
\frac{\k[x_0,x_1,y_2\ldots,y_{g}]}{(y_2-f_2(x_0,x_1),\ldots,y_g-f_g(x_0,x_1,y_2,
\ldots,y_{g-1}))}\simeq \k[x_0,x_1].$$
Let $\Sigma'$ be a regular subdivision of $\mathbb{R}_+^{g+1},$ which is 
compatible with the Newton dual fan of $y_i-f_i,,i=2,\ldots,g$ (see section 5 of 
\cite{GT} for the construction of $\Sigma'$), and let $\Sigma''$ be the Stellar 
subdivision of $\Sigma'$ associated with the vector  $v_{\nu_E}.$ Finally let 
$\Sigma$ be a regular subdivision of $\Sigma''.$ Then the first 3 properties of 
the proposition follows from theorem 5.2 in \cite{GT}. Now by construction of 
the embedding $e,$ we have that if $\mathbb{W}$ is the irreducible component of 
$Cont^p(\mathcal{C})$ which defines $\nu_E,$ then
$$e_{\infty}(\mathbb{W})=e_{\infty}(\A^2_{\infty})\cap 
Cont^{\nu_E(x_0)}(x_0)\cap Cont^{\nu_E(x_1)}(x_1)\cap 
Cont^{\nu_E(x_2)}(y_2)\cap\ldots\cap Cont^{\nu_E(x_g)}(y_g),$$

where $e_{\infty}:\A^2_{\infty}\hookrightarrow \A^{g+1}_{\infty}$ is the 
canonical morphism. But the divisorial valuation associated with 
$$\mathbb{U}=Cont^{\nu_E(x_0)}(x_0)\cap Cont^{\nu_E(x_1)}(x_1)\cap 
Cont^{\nu_E(x_2)}(y_2)\cap\ldots\cap Cont^{\nu_E(x_g)}(y_g)\subset 
\A^{g+1}_{\infty}$$ is $\nu_{E'},$ which in terms of arcs means that  
$\mu_{\infty}(Cont^1(E'))$ dominates $\mathbb{U},$ hence we 
have that $\eta_{\infty}(Cont^1(E))$ dominates $\mathbb{W}$ where $\eta$ is the 
restriction of $\mu$ to  $\tilde{\A^2}.$ The property 4 in the proposition 
follows form the description of the valuation associated to $\mathbb{W}.$

\end{proof}

\begin{rem} Note that that we can use the equations $f_i$ to define an 
overweight deformation in the sense of $\cite{T2},$ hence $\nu_E$ can 
be obtained from the monomial valuation $\nu_{E'}$ as in  proposition 3.3 
\cite{T2}.  
\end{rem}
\begin{ex}
Let $\mathcal{C}$ be the irreducible curve defined by the equation 
$x_1^2-x_0^3=0.$ Let $\nu$ be the valuation defined by $C_6\subset 
\mathcal{C}_6$ or equivalently by the corresponding irreducible component of 
$Cont^7(x_1^2-x_0^3).$  Note that the ideal of $C_6$ is generated by 
$$(x_0^{(0)},x_0^{(1)},x_1^{(0)},\ldots,x_1^{(2)},{x_1^{(3)}}^2-{x_0^{(2)}}
^3).$$ Then by the discussion at the beginning of this section, we have that 
$x_0,x_1$ and $x_2=x_1^2-x_0^3$ give a minimal generating sequence of $\nu.$  We 
embed $\A^2=\mbox{Spec} \k[x_0,x_1]\hookrightarrow \A^3=\mbox{Spec} 
\k[x_0,x_1,y_2]$ by the equation  $y_2-(x_1^2-x_0^3)=0.$ A subdivision of 
$\mathbb{R}_+^3$ as in the proposition \ref{emb} is given by a fan $\Sigma$ 
whose edge vectors are the vectors $$(1,1,1),(1,2,3),(2,3,5),(2,3,6),(2,3,7)$$ 
where the last vector is the $v_{\nu}=(\nu(x_0),\nu(x_1),\nu(x_2)).$ We are 
interested in a chart of $X_\Sigma$ where we can see the divisor $E'$ associated 
with the vector  $v_{\nu}.$ We consider the chart $X_\sigma = \A^3=\mbox{Spec} 
\k[u,v,w]$ generated by the vectors $(1,2,3),(2,3,6),(2,3,7).$ The restriction 
of $\mu$ to this chart is given by $$ x_0=uv^2w^2 $$
 $$x_1=u^2v^3w^3 $$
  $$y_2=u^3v^6w^7.$$
The strict transform of $\A^2=\{y_2-(x_1^2-x_0^3)=0\}\subset \A^3$ is given by 
$$\tilde{\A^2}=\{w-u+1=0\}\simeq \mbox{Spec} \k[u,v]\subset \A^3=\mbox{Spec} 
\k[u,v,w]$$
and $E'$ is defined by $w=0.$ Thus the divisor $E$ is defined in $\tilde{\A^2}$ 
by the equation $u-1=0.$ The restriction $\eta$ of $\mu$ to $\tilde{\A^2}$ is 
obtained from the description of $\mu$ by substituting $w$ by $u-1.$ Hence 
$\eta$ is given by 

$$ x_0=uv^2(u-1)^2 $$
 $$x_1=u^2v^3(u-1)^3.$$
 It is direct to verify that $\eta$ is obtained as follows: First we we consider 
the minimal embedded resolution of the curve $\mathcal{C}=\{x_1^2-x_0^3=0 \}$ 
(which is obtained by three consecutive point blowing ups), then we blow up  the 
intersection of the  strict transform of $\mathcal{C}$ with the exceptional 
divisor. The divisor obtained from this last blowing up satisfies $\nu=\nu_E.$ 
We see that the total transform of $\mathcal{C}$ by $\eta$ is given by the 
equation $u^3v^6(u-1)^7$ and hence  that $\nu_E(x_1^2-x_0^3)=7.$
 \end{ex}

This result shows a different approach from the one of \cite{GT} to the 
resolution of singularities of an irreducible plane curve $\mathcal{C}$ by one 
toric morphism. Indeed in loc.cit. the embedding $e$ is constructed  from the 
study of the curve valuation $\nu_{\mathcal{C}},$ while the approach suggested 
by this article is to study the divisorial valuation associated with the 
irreducible component  $C_{p-1}$ of $\mathcal{C}_{p-1}$(where 
$p=n_g\bar{\beta}_g$ is detected via invariants of jet schemes). The two 
approaches lead to the same embedding in this case, in higher dimensions they 
bifurcate.\\ 

Let us explain a little bit more the point of view suggested in this article 
about the embedding $e.$  Let $\nu=\nu_{\alpha}$ be the monomial valuation 
defined on $\A^n=Spec \k[x_1,\ldots,x_n] $ by a vector 
$\alpha=(\alpha_1,\ldots,\alpha_n),$ where $\alpha_i \in 
\mathbb{N},i=1,\ldots,n.$ Let  $I\subset \k[x_1,\ldots,x_n]$ be an ideal and we 
assume that the origin $O$ belongs to the variety $V(I)\subset \A^n=Spec 
\k[x_1,\ldots,x_n]$ defined by this ideal. We will say that $I$ or $V(I)$ is 
non-degenerate with respect to $\nu$ at $O$ if the singular locus of the variety 
defined by the initial ideal $in_{\nu}(I)$ of $I$ does not intersect the torus 
$(\k^*)^n.$ Note that in this context, the initial ideal of $I$ is defined by
$$ in_{\nu}(I)=\{in_{\nu}(f),f\in I \},$$
where for  $f=\sum a_{i_1,\ldots,i_n}x_1^{i_1}\cdots x_n^{i_n}\in 
\k[x_1,\ldots,x_n],$ 
$$in_{\nu}(f)=\sum_{a_{i_1,\ldots,i_n}\not=0, 
i_1\alpha_1+\cdots+i_n\alpha_n=\nu(f)} a_{i_1,\ldots,i_n}x_1^{i_1}\cdots 
x_n^{i_n}.$$

It follows from \cite{AGS},\cite{Te1} (see also \cite{Va} for the hypersurface 
case) that if for every $\alpha=(\alpha_1,\ldots,\alpha_n),\alpha_i \in 
\mathbb{N},i=1,\ldots,n,$ $I$ is non-degenerate with respect to $\nu_{\alpha}$ 
at $O,$ then we can construct a proper toric birational morphsim 
$Z\longrightarrow \A^n$ that resolve the singularities of $V(I)$ in a 
neighborhood of $O.$ Notice that $I$ can be degenerate with respect to a 
valuation defined by a vector $\alpha$ if there exists an irreducible family of 
jets (having a large contact with $V(I)$) or arcs on $V(I)$ which satisfy that  
for a generic $\gamma=(\gamma_1(t),\ldots,\gamma_n(t))$ in this family, $(ord_t 
\gamma_1(t),\ldots,ord_t\gamma_n(t))=\alpha:$ indeed, by Newton-Puiseux type 
theorem, if this is not satisfied,   $in_{\nu_{\alpha}}(f)$ will contain 
monomials, hence by definition $I$ will be non-degenerate with respect to 
$\nu_{\alpha}.$ By studying irreducible components of jet schemes of a plane 
branch $\mathcal{C},$ as we have done, we are also looking for the degeneracy 
behind the first Newton polygon.  The embedding we have constructed by applying 
proposition \ref{emb} to the divisorial valuation associated with the 
irreducible component $C_{n_g\bar{\beta}_g-1}$ of 
$\mathcal{C}_{n_g\bar{\beta}_g-1},$ have the following  property: Let $I$ be the 
defining ideal of the curve $\mathcal{C}$ in $\A^{g+1}$ and let 
$\alpha=(\bar{\beta}_0,\ldots,\bar{\beta}_g);$  the initial ideal  
$in_{\nu_{\alpha}}$(I) is the defining ideal of the monomial curve defined by 
$\{ (t^{\bar{\beta}_0},\ldots,t^{\bar{\beta}_g}), t \in \k \},$ which has an 
isolated singularity at $O,$
hence $I$ is non-degenerate with respect to  $\nu_{\alpha}.$ Moreover, this is 
the only relevant vector $\alpha$ with respect to which we should check 
degeneracy, the reason being that the initial ideal with respect to any other 
vector will contain monomials. One crucial thing is that in the curve case, the 
initial ideal we found is binomial, thus it defines a toric variety, in higher 
dimension it will not be the case, and more technology will be needed.

Hussein Mourtada,\\
Equipe G\'eom\'etrie et Dynamique, \\
Institut Math\'ematique de Jussieu-Paris Rive Gauche,\\
 Universit\'e Paris 7, \\
 B\^atiment Sophie Germain, case 7012,\\
75205 Paris Cedex 13, France.\\

Email: hussein.mourtada@imj-prg.fr
 
\end{document}